\documentclass [12pt] {amsart}
\usepackage{amsmath, amsfonts, latexsym, array, epsfig, amssymb, psfrag}
\usepackage[all]{xypic}
\usepackage{graphicx}
\newtheorem{thm}{Theorem}[section]
\newtheorem{cor}[thm]{Corollary}
\newtheorem{lem}[thm]{Lemma}

\newtheorem{que}[thm]{Question}

\begin{document}
\title[Centralizer]
{Trivial centralizers for codimension-one attractors}
\author{Todd Fisher}
\address{Department of Mathematics, Brigham Young University, Provo, UT 84602}
\email{tfisher@math.byu.edu}
\thanks{}

\date{November 2, 2007}
\commby{}

\begin{abstract}
 We show that if $\Lambda$ is a codimension-one hyperbolic attractor for a $C^r$ diffeomorphism $f$, where $2\leq r\leq \infty$, and $f$ is not Anosov, then there is a neighborhood $\mathcal{U}$ of $f$ in $\mathrm{Diff}^r(M)$ and an open and dense set $\mathcal{V}$ of $\mathcal{U}$ such that any $g\in\mathcal{V}$ has a trivial centralizer on the basin of attraction for $\Lambda$.    \end{abstract}

\maketitle

\section{Introduction}

Let $f\in\mathrm{Diff}^r(M)$ be a $C^r$ diffeomorphism where $M$ is a connected compact manifold  without boundary.  The {\it centralizer} of $f$, denoted $Z(f)$, is the set of all diffeomorphisms that commute with $f$: 
 $$Z(f)=\{g\in\mathrm{Diff}^r(M)\,|\, fg=gf\}.$$ Note that the set $Z(f)$ contains  all powers of $f$.   
 The set $Z(f)$ forms a group under composition.   We say $f$ has {\it trivial centralizer} if  $Z(f)=\{f^n\, |\, n\in\mathbb{Z}\}$.   A diffeomorphism with trivial centralizer has no smooth symmetries.  Smale~\cite{Sma2, Sma} asked the following question concerning the frequency of diffeomorphisms with trivial centralizer.

\begin{que}\label{q.smale}
Let $1\leq r\leq \infty$, $M$ be a smooth, connected, compact, boundaryless manifold, and $T$ be the set of all $f\in\mathrm{Diff}^r(M)$ such that $f$ has trivial centralizer.
\begin{enumerate}
\item  Is $T$ dense in $\mathrm{Diff}^r(M)$?
\item Is $T$ residual in $\mathrm{Diff}^r(M)$?
\item Does $T$ contain an open and dense subset of $\mathrm{Diff}^r(M)$?
\end{enumerate}
\end{que}

There are a number of partial results to the above question.  A few of the results are the following:
\begin{itemize}
\item   Kopell~\cite{Kop1} shows that $T$ contains an open and dense set if  $M=S^1$ and $r=\infty$.  \item Palis and Yoccoz~\cite{PY1} show there is an open and dense set of $C^{\infty}$ Axiom A diffeomorphisms with the strong transversality property and a periodic sink that have a trivial centralizer.  
\item In~\cite{Fis2} it is shown that there is a residual set of $C^\infty$ Axiom A diffeomorphisms with the no cycles property with trivial centralizer.  Additionally, for a surface and $2\leq r\leq\infty$  it is shown that there exists an open and dense set of $C^r$ Axiom A diffeomorphisms with the no cycles property whose elements have trivial centralizer.  
\item Burslem~\cite{Bur} shows that for certain classes of $C^{\infty}$ partially hyperbolic diffeomorphisms there is a residual set with trivial centralizer.  
\item Bonatti, Crovisier, and Wilkinson~\cite{BCW1} have shown that $C^1$ diffeomorphisms contain a residual set with trivial centralizer.
\item Bonatti, Crovisier, Vago, and Wilkinson~\cite{BCW} are able to show that on any manifold there is a open set in $\mathrm{Diff}^1(M)$ containing a dense set whose elements do not have a trivial centralizer.
\end{itemize}

One approach to Question~\ref{q.smale} is to examine the centralizer of a diffeomorphism restricted to the basin of attraction for an attractor.  This is the approach used in~\cite{PY1} and ~\cite{Fis2}.  The next theorem addresses this approach in the case of codimension-one attractors.

\begin{thm}\label{t.codimone}
If $\Lambda$ is a codimension-one hyperbolic attractor for a $C^r$
diffeomorphism $f$ where $2\leq r\leq \infty$ and $f$ is not Anosov, then there exists a neighborhood $\mathcal{U}$ of $f$ in
$\mathrm{Diff}^r(M)$ and an open and dense set $\mathcal{V}$ of $\mathcal{U}$ such that for all $g\in \mathcal{V}$ and $h\in Z(g)$
the map $h|_{W^s(\Lambda(g))}=g^n|_{W^s(\Lambda(g))}$ where $\Lambda(g)$ is the
continuation of $g$.
\end{thm}

We note that the novelty of the above theorem is that $r$ can be taken less than infinity and the the fact that $\mathcal{V}$ is open.  The above result will also hold for codimension-one hyperbolic repellers.

Two cases where we can apply the above result to the entire manifold are for certain structurally stable diffeomorphisms and quasi-Anosov diffeomorphisms of $3$-manifolds.

A function $f$ is {\it structurally stable} if there exists a neighborhood $\mathcal{U}$ of $f$ such that for each $g\in \mathcal{U}$ there exists a topological conjugacy from $f$ to $g$.  From the work of Grines and Zhuzhoma~\cite{GZ} on non-Anosov structurally stable diffeomorphisms with a codimension-one attractor we know in this case that the basin of attraction of the hyperbolic attractor is dense in the manifold.  We then have the following corollary.

\begin{cor} If $f$ is a structurally stable, non-Anosov, $C^r$ diffeomorphism where $2\le r\leq \infty$, and contains a codimension-one attractor, then there exists a neighborhood $\mathcal{U}$ of $f$ in $\mathrm{Diff}^r(M)$ and an open and dense set $\mathcal{V}$ of $\mathcal{U}$ whose elements have a trivial centralizer.
\end{cor} 

A diffeomorphism $f$  is {\it expansive} if there exists a constant $c>0$  such that for any two distinct points $x$ and $y$ there exists an integer $n$ where $d(f^n (x),f^n(y))\geq c$.  The definition is independent of the metric, although the constant $c$ in general depends on the metric.
A diffeomorphism $f$ is {\it quasi-Anosov} if there exists a neighborhood $\mathcal{U}$ of $f$ such that each $g\in \mathcal{U}$ is expansive.   
For quasi-Anosov diffeomorphisms we will show the next corollary follows from results in~\cite{FRH} and Theorem~\ref{t.codimone}.

\begin{cor} \label{c.2} If $\mathrm{dim}(M)=3$, then there is a $C^r$ open and dense set of $C^r$ non-Anosov quasi-Anosov diffeomorphisms of $M$ where $2\leq r\leq \infty$ whose elements have a trivial centralizer.
\end{cor}

From the above results the next question is a natural weakening of   Question~\ref{q.smale} part (3). 
\begin{que}
If $2\leq r\leq \infty$, then does $T$ contain an open and dense subset of $\mathrm{Diff}^r(M)$?
\end{que}

\section{Background}

We now review some basic definitions and facts about hyperbolic attractors and commuting diffeomorphisms.
We assume that all of our maps are diffeomorphisms of a manifold to itself.

A compact set $\Lambda$ invariant under the action of $f$ is {\it
hyperbolic} if there exists a splitting of the  tangent space
$T_{\Lambda}f=\mathbb{E}^u\oplus \mathbb{E}^s$ and positive
constants $C$ and $\lambda<1$ such that, for any point
$x\in\Lambda$ and any $n\in\mathbb{N}$,
$$
\begin{array}{llll}
\| Df_{x}^{n}v\|\leq C \lambda^{n}\| v\|,\textrm{ for }v\in
E^{s}_x \textrm{, and}\\
\| Df_{x}^{-n}v\|\leq C \lambda^{n}\| v\|,\textrm{ for }v\in
E^{u}_x. \end{array}
$$

For $\epsilon>0$ sufficiently small and $x\in \Lambda$ the
\textit{local stable and unstable manifolds} are respectively:
$$
\begin{array}{llll}
W_{\epsilon}^{s}(x,f)=\{ y\in M\, |\textrm{ for all }
n\in\mathbb{N}, d(f^{n}(x), f^{n}(y))\leq\epsilon\},\textrm{
and}\\
W_{\epsilon}^{u}(x,f)=\{ y\in M\, |\textrm{ for all }
n\in\mathbb{N}, d(f^{-n}(x), f^{-n}(y))\leq\epsilon\}.
\end{array}$$
The \textit{stable and unstable manifolds} are respectively:
$$
\begin{array}{llll}
W^s(x,f)=\bigcup_{n\geq 0}f^{-n}\left(
W_{\epsilon}^s(f^n(x),f)\right) ,\textrm{ and}\\
W^u(x,f)=\bigcup_{n\geq
0}f^{n}\left(W_{\epsilon}^u(f^{-n}(x),f)\right). \end{array}
$$
For a $C^r$ diffeomorphism the stable and unstable manifolds of a hyperbolic set are $C^r$
injectively immersed submanifolds.

A useful property is the following standard result.

\begin{thm}(Structural stability of hyperbolic sets) Let $f\in \mathrm{Diff}(M)$ and $\Lambda$ be a hyperbolic set for $f$.  Then for any neighborhood $V$ of $\Lambda$ and every $\delta>0$ there exists a neighborhood $\mathcal{U}$ of $f$ in $\mathrm{Diff}(M)$ such that for any $g\in \mathcal{U}$ there is a hyperbolic set $\Lambda_g\subset V$ and a homeomorphism $h:\Lambda_g\rightarrow \Lambda$ with $d_{C^0}(\mathrm{id},h) + d_{C^0}(\mathrm{id},h^{-1})<\delta$ and $h\circ g|_{\Lambda_g}=f|_{\Lambda}\circ h$.  Moreover, $h$ is unique when $\delta$ is sufficiently small.
\end{thm}

Many of the results in this paper will follow from looking at the basins of attractors and repellers.  We now give some basic properties of these sets.

A set $X\subset M$ has an {\it attracting neighborhood} if there
exists a neighborhood $V$ of $X$ such that
$X=\bigcap_{n\in\mathbb{N}}f^n(V)$.  
If $X$ is a compact
set of a smooth manifold $M$ and $f$ is a continuous map from $M$
to itself, then $f|_X$ is {\it transitive} if for any open
sets $U$ and $V$ of $X$ there exists some $n\in\mathbb{N}$ such
that $f^n(U)\cap V\neq\emptyset$. 
A set $\Lambda\subset M$
is called a {\it hyperbolic attractor} 
if $\Lambda$ is a transitive hyperbolic set for a diffeomorphism
$f$ with an attracting neighborhood.  The {\it basin of attraction} for a hyperbolic attractor $\Lambda$ is the set 
$$W^s(\Lambda)=\{x\in M\, |\, d(f^n(x), \Lambda)\rightarrow 0\textrm{ as }n\rightarrow\infty\}.$$
A standard result is that if $p$ is a periodic point contained in a hyperbolic attractor $\Lambda$, then
$$W^s(\Lambda)\subset \overline{W^s(\mathcal{O}(p))}.$$

We now review some basic properties of commuting diffeomorphisms.
Let $f$ and $g$ be commuting diffeomorphisms.  Let
$\mathrm{Per}^n(f)$ be the periodic points of period $n$ for $f$
and $\mathrm{Per}^n_h(f)$ denote the hyperbolic periodic points in
$\mathrm{Per}^n(f)$.  If $p\in\mathrm{Per}^n(f)$, then
$g(p)\in\mathrm{Per}^n(f)$ so $g$ permutes the points of
$\mathrm{Per}^n(f)$.  Furthermore, if
$p\in\mathrm{Per}^n(f)$, then
$$T_{g(p)}f^nT_pg=T_pgT_pf^n.$$
Hence, the linear maps $T_{g(p)}f^n$ and $T_pf^n$ are similar.
Since $\#(\mathrm{Per}^n_h(f))<\infty$ it
follows that if $p\in\mathrm{Per}^n_h(f)$, then $g(p)\in\mathrm{Per}(g)$.  If $p\in\mathrm{Per}^n_h(f)$, then
$$g(W^u(p,f))=W^u(g(p),f)\textrm{ and }g(W^s(p,f))=W^s(g(p),f).$$

\section{Trivial centralizer for codimension-one attractors}

Let $p$ denote a periodic point in a codimension-one attractor.  The proof of Theorem~\ref{t.codimone} will depend on a linearization of $W^s(p)$ to the reals.  Since $f$ is in $\mathrm{Diff}^r(M)$ where $2\leq r\leq\infty$ we know there exists a $C^2$ immersion $\phi$ of $W^s(p)$ such that $F=\phi f\phi^{-1}\in \mathrm{Diff}^r(\mathbb{R})$ fixing the origin.  The following result due to Sternberg will then help $C^1$ linearize the map $F$.

\begin{thm}~\cite{Ste} Let $f$ be a $C^r$ monotone increasing function defined in a neighborhood of the origin in $\mathbb{R}$ where $2\leq r\leq\infty$ and $|f'(0)|=a\neq 0,1$.  Then there exists a function $g$ that is a $C^{r-1}$ monotone increasing function with $C^{r-1}$ inverse defined in a neighborhood of the origin in $\mathbb{R}$ such that $g^{-1} f g(x)=ax$ for all sufficiently small $x$ in $\mathbb{R}$.
\end{thm}

As the function $F$ above is defined on all of $\mathbb{R}$ and globally contracting we will want the following corollary to the result of Sternberg.

\begin{cor}\label{c.linear} Let $f\in\mathrm{Diff}^r(\mathbb{R})$ be a contraction mapping fixing the origin where $2\leq r\leq\infty$ and $|f'(0)|=a\neq 0,1$.  Then there exists a function $g$ in $\mathrm{Diff}^{r-1}(\mathbb{R})$ such that $g^{-1} f g(x)=ax$ for all  $x$ in $\mathbb{R}$.
\end{cor}

For the proof of the corollary one simply takes the $C^{r-1}$ linearizing function defined in a neighborhood of the origin given by Sternberg and extends to all of the reals through iteration.  Since
the map is a contraction mapping this can be extended to all of the reals.\\





\noindent{\bf Proof of Theorem~\ref{t.codimone}.}
Let $f$ be a diffeomorphism and $\Lambda$ a codimension-one hyperbolic attractor.  By the structural  stability of hyperbolic sets there exists an open set $\mathcal{U}$ of $f$ in $\mathrm{Diff}^r(M)$ such that for all $g\in \mathcal{U}$ there is a continuation $\Lambda(g)$ of $\Lambda$.  Since $\Lambda$ is a transitive locally maximal hyperbolic set the periodic points are dense, see~\cite[p. 574]{KH1}.

From~\cite{GZ} there exists a periodic point $p\in\Lambda$, called a periodic boundary point such that $W^s(p)\cap\Lambda$ is a Cantor set in one of the separtrices and the empty set in the other.   Furthermore, the map $f|_{W^s(p)\cap \Lambda}$ is invariant so $f|_{W^s(p)}$ is orientation preserving.

For simplicity we assume that $p$ is a fixed point, we leave the details of $p$ a periodic point to the reader.
We let $\mathcal{V}$ be an open and dense set  in $\mathcal{U}$ such that if $g\in\mathcal{V}$ and $q$ and $q'$ are hyperbolic fixed points of $g$, then $T_{q}g$ is not similar to $T_{q'}g$.  Hence, if $g\in \mathcal{V}$ and $h\in Z(g)$, then $h(p(g))=p(g)$.  

The proof of Theorem~\ref{t.codimone} will follow from examining $W^s(p)$ since this is dense in $W^s(\Lambda)$.
The set $W^s(p)$ is a $C^2$ immersed copy of $\mathbb{R}$ where the map $\psi:W^s(p)\rightarrow\mathbb{R}$ is the immersion with $\psi(p)=0$.  Then the map $G=\psi g\psi^{-1}$ is a $C^2$ contracting diffeomorphism of the reals fixing the origin and $G'(0)=a\in(0,1)$ (since $G$ is orientation an orientation preserving contraction).
From Corollary~\ref{c.linear} there exists a map $\phi\in\mathrm{Diff}^1(M)$ such that $(\phi G\phi^{-1})(x)=ax$.

The next lemma will show $h|_{W^s(p)}$  is also $C^1$ linearized by $\Phi=\phi\psi$.

\begin{lem}\label{l.commutereals}~\cite[p. 61]{KH1}  Any $C^1$ map defined on a neighborhood of the origin on the real line and commuting with a linear contraction is linear.
\end{lem}

Let $L=\Phi g\Phi^{-1}$ and $H=\Phi h\Phi^{-1}$.  We now show that $LH=HL$ for  all $x$.  
$$
\xymatrix{
\mathbb{R}\ar[d]_{\Phi^{-1}}\ar[r]^{L} &  \mathbb{R}\ar[d]^{\Phi^{-1}}\\
W^s(p(g))\ar[d]_{h}\ar[r]^{g} & W^s(p(g))\ar[d]^{h}\\
W^s (p(g))\ar[d]_{\Phi}\ar[r]^{g} & W^s (p(g))\ar[d]^{\Phi}\\
\mathbb{R}\ar[r]_{L} & \mathbb{R}}
$$
Since $L$ is linear and $H$ is $C^1$ we know from the previous lemma that $H$ is linear.

To complete the proof of Theorem~\ref{t.codimone} we show there exists an open and dense set $\mathcal{V}'$ in $\mathcal{V}$ such that for each $g\in \mathcal{V}'$ and $h\in Z(g)$ we have $h|_{W^s(p(g))}=g^n|_{W^s(p(g))}$ for some $n$.  
Since $W^s(p(g))$ is dense in $W^s(\Lambda)$ we see that $h|_{W^s(\Lambda)}=g^n|_{W^s(\Lambda)}$.  To construct the set $\mathcal{V}'$  we perturb the Cantor sets on one separatrix of $W^s(p(g))$ so that the only linear maps that commute with $L$ are powers of $L$.

Define $C=\Phi(W^s(p)\cap\Lambda)$.  So $C$ is a Cantor set in the reals.  We may assume that $C\subset[0,\infty)$.   A {\it gap} of $C$ is a connected component of $\mathbb{R}-C$.  A bounded gap is a bounded component.  

\begin{figure}[htb]
\begin{center}
\psfrag{x}{$x$}
\psfrag{y}{$y$}
\psfrag{l}{$L(U)$}
\psfrag{z}{$L(x)$}
\psfrag{w}{$L(y)$}
\psfrag{1}{$U_1$}
\psfrag{2}{$U_2$}
\psfrag{n}{$U_n$}
\psfrag{u}{$U$}

\includegraphics{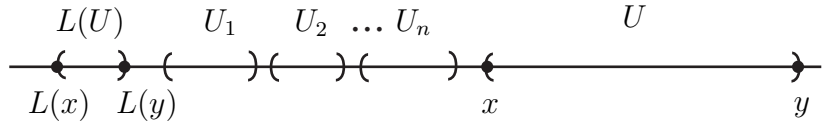}
\caption{Gaps between $U$ and $L(U)$}\label{f.gaps1}
\end{center}
\end{figure}

Let $U$ be a bounded gap in $C$ with boundary points $0<x<y<\infty$.    Let $U_1,...U_n$ be the gaps of size at least size
$L(U)$ contained in the interval $[L(y),x]$ and $x_i<y_i$ be the boundary points of gap $U_i$ for $1\leq
i\leq n$ where $x_i<x_{i+1}$ for all $1\leq i<n$.  We know there exists some $k\in\mathbb{Z}$ such that $L^kH(x)\in [L(x), x)$.  Furthermore, we know that $C$ is $H$ invariant so $L^kH(x)\in C$ and $L^kH(U)$ is a gap of $C$ contained in $[L(x), x)$.  Since $L^kH$ is a linear map we know the size of $L^kH(U)$  by the image of $x$ under $L^kH$.  Additionally, $L^kH(U)$ is  contained in $[L(x), x)$ and the linearization implies that the size of $L^kH(U)$ is at least that of $L(U)$.  Hence, $L^kH(U)$ is either $L(U)$ or one of the set $U_1,...,U_n$. 

If $L^kH(U)=U_i$, then $x_i/x=y_i/y$ is the contraction constant of $L^kH$.  By an arbitrarily small perturbation of $g$ we can ensure that $x_i/x\neq y_i/y$ for all $1\leq i\leq n$.  
Indeed, let $W=\Phi^{-1}([0, y])\subset W^s(p(g))$ and $V$ be a neighborhood of $W$ and $\epsilon>0$ such that $$V\cap \bigcup_{i=1}^nB_{\epsilon}(g^{-1}\Phi^{-1}(y_i))=\emptyset.$$  Then there exists an arbitrarily small perturbation of $g$ with support in $\bigcup_{i=1}^nB_{\epsilon}(g^{-1}\Phi^{-1}(y_i))$ that does not change the linearization of $W$ and ensures that $x_i/x\neq y_i/y$ for all $1\leq i\leq n$.

Furthermore, this will be an open condition since $L$ varies continuously with $g$.
Hence, $L^kH(U)=L(U)$ and $H$ is a power of $L$ showing that $h|_{W^s(\Lambda)}=g^{1-k}|_{W^s(\Lambda)}$.
$\Box$\\


\noindent{\bf Proof of Corollary~\ref{c.2}.}  From~\cite{FRH} we know that for every quasi-Anosov diffeomorphism of a 3-manifold there is an open and dense set of points in the manifold contained in the basin of a codimension-one attractor or repeller.  

Let $QA^r(M)$ be the set of quasi-Anosov diffeomorphisms on a manifold $M$ that are non-Anosov.
Since each quasi-Anosov diffeomorphism is Axiom A with no cycles~\cite{Man1} we know from Proposition 3.2 in~\cite{Fis2} that there is an open and dense set $\mathcal{U}$ in $QA^r(M)$ such that for each  $f\in\mathcal{U}$ and $\Lambda$ and $\Lambda'$ are attractors for $f$ where $$\overline{W^s(\Lambda)}\cap \overline{W^s(\Lambda')}\neq\emptyset,$$ then there exists a hyperbolic repeller $\Lambda_r$ such that $$W^s(\Lambda)\cap W^u(\Lambda_r)\neq\emptyset\textrm{ and }W^s(\Lambda')\cap W^u(\Lambda_r)\neq\emptyset.$$

From Theorem~\ref{t.codimone} there exists an open and dense set $\mathcal{V}\subset \mathcal{U}$ such that if $f\in\mathcal{V}$, $\Lambda$ a codimension-one hyperbolic attractor or repeller of $f$, and $g\in Z(f)$, then $g|_{W^s(\Lambda)}=f^n|_{W^s(\Lambda)}$ if $\Lambda$ is an attractor or $g|_{W^u(\Lambda)}=f^n|_{W^u(\Lambda)}$ if $\Lambda$ is a repeller where $n\in\mathbb{Z}$.

Let $f\in\mathcal{V}$, $\Lambda_a$ be a codimension-one hyperbolic attractor for $f$, $\Lambda_r$ be a codimension-one repeller for $f$, $W^s(\Lambda_a)\cap W^u(\Lambda_r)\neq \emptyset$, and $g\in Z(f)$.  Let $$g|_{W^s(\Lambda_a)}=f^n|_{W^s(\Lambda_a)}\textrm{ and  }g|_{W^u(\Lambda_r)}=f^m|_{W^u(\Lambda_r)}.$$  Suppose that  $m\neq n$.  Then 
$$g'=f^{-n}g\in Z(f)\textrm{ and }g'|_{W^s(\Lambda_a)}=\mathrm{id}|_{W^s(\Lambda_a)}.$$  
Let $x\in W^s(\Lambda_a)\cap W^u(\Lambda_r)$.  Then 
$$g'(x)=x=f^{m-n}(x).$$  Hence, $x$ is periodic for $f$, a contradiction.  Therefore, we know for each $f\in \mathcal{V}$ there is an open and dense set $U\subset M$ such that if $g\in Z(f)$, then $g|_U=f^n|_U$.  This follows from Proposition 3.2 in~\cite{Fis2}.  Hence, for each $f\in\mathcal{V}$ and $g\in Z(f)$ we have $g=f^n$. $\Box$

\bibliographystyle{plain}
\bibliography{commute012808}

\end{document}